\documentclass{amsart}

\usepackage{amssymb}
\usepackage{amsmath}
\usepackage{url}
\usepackage{xcolor}
\usepackage{booktabs}
\usepackage[english]{babel}
\usepackage[linesnumbered,norelsize,ruled,vlined]{algorithm2e}
\usepackage{algorithmic}
\usepackage{multirow}
\usepackage[quotecolor=red]{revquote}
\usepackage[normalem]{ulem}
\usepackage{graphicx}
\usepackage[latin1]{inputenc}
\usepackage{natbib}
\bibpunct{(}{)}{;}{a}{,}{;}

\begin{document}

\title[Induced p-median problem with upgrades]{Metaheuristic algorithms 
for the induced \textit{P}-median problem with upgrades}

\author[S. Salazar]{Sergio Salazar}
\address[S. Salazar]{
Universidad Rey Juan Carlos, Mostoles, Madrid, Spain
}
\email[S. Salazar] {sergio.salazar@urjc.es}

\author[A. Duarte]{Abraham Duarte}
\address[A. Duarte]{
Universidad Rey Juan Carlos, Mostoles, Madrid, Spain
}
\email[A. Duarte]{abraham.duarte@urjc.es}

\author[M.G.C. Resende]{Mauricio G.C. Resende}
\address[M.G.C. Resende]{ITA-UNIFESP Graduate Program in Operations Research,
            S\~ao Jos\'e dos Campos, S\~ao Paulo, Brazil and
Center for Discrete Mathematics and Theoretical
Computer Science (DIMACS),
            Piscataway,
            New Jersey, USA
}
\email[M.G.C. Resende]{mauricio.resende@gp.ita.br, mgcr@dimacs.rutgers.edu}

\author[J.M. Colmenar]{J. Manuel Colmenar}
\address[J.M. Colmenar]{
Universidad Rey Juan Carlos, Mostoles, Madrid, Spain
}
\email[J.M. Colmenar]{josemanuel.colmenar@urjc.es}

\begin{abstract}
Facility location problems (FLPs) are a family of optimisation
problems with significant social impact. This class of problems has
been the subject of study since the 1960s, with classical approaches
including the Weber problem and the \textit{p}-Median problem.
Currently, more complex variations of these problems are being
investigated. In particular, the Induced \textit{p}-Median Problem
with Upgrades (IpMU) represents a variation of the classical
\textit{p}-Median problem, where the concepts of transport cost and
time are separated as distinct metrics in the input graph of the
problem. Furthermore, the problem includes a budget which allows
one to relax the graph costs, reducing the cost of the edges, thus
improving the associated routes between the designated medians and
the customers. In this study, a metaheuristic algorithm, based on
the Greedy Randomized Adaptive Search Procedure (GRASP), is proposed.
A two-phase resolution scheme is defined, studying the median problem
and the upgrading problem independently. In this approach, a larger
set of state-of-the-art instances was analysed to ensure a fair
comparison with previous proposals. In addition, the characteristics
of the instances were studied to assess their complexity. The results
obtained are promising when compared to the state-of-the-art, which
is based entirely on mathematical programming models. The execution
time was improved on average by two orders of magnitude for the
harder instances, and the best known result was obtained in more
than 99\% of the tested instances.
\end{abstract}

\date{February 17, 2026}
\thanks{Technical Report. Paper to appear in 
\textit{Knowledge-Based Systems}.}
\maketitle

\section{Introduction}
\label{sec:intro}
% Revisada
One of the most studied topics in Operations Research is the family
of Facility Location Problems (FLPs)
\citep{francis1992facility,drezner2004facility}. The interest in
this family of problems is twofold. On the one hand, their mathematical
analysis has promoted different solution methods, from those based
on mathematical models to those based on advanced metaheuristic
algorithms \citep{drezner2004facility,laporte2020introduction}. On
the other hand, the impact of this family of problems in real-life
scenarios is quite straightforward. The location of ammunition
depots at the end of armed conflicts to avoid the risks associated
with their life cycle \citep{daugistanli2024facility} and the location
and distribution of medical facilities during the COVID-19 pandemic
\citep{fan2022distributionally} are two recent examples of this
variety of applied scenarios.

% \sout{
% The FLP family is extensive and its taxonomy is based on a variety of different characteristics. Taking into account the mathematical characteristics of the problem model, unilevel and multilevel problems (MFLPs) can be found in the literature \cite{ortiz2018multi}, where several objective functions correspond to different constraints in the complete model. Considering the mathematical characteristics, problems are classified as discrete or continuous, depending on whether the domain of solutions, and therefore the search space, is discrete or continuous. In the case of combinatorial optimization, this family includes classical problems such as the \textit{p}-Median problem \cite{alp2003efficient} and the \textit{p}-Center problem \cite{plesnik1987heuristic}. In continuous space, other problems can be found, such as the Fermat-Weber problem \cite{cooper1968extension}.}

Within the field of FLP, we can find different families of problems
which present specialized characteristics, objectives, and constraints.
In this regard, Healthcare Facility Location Problems (HCFLPs),
related to disaster management and humanitarian logistics, emerges
as one of the most studied families. Among them, a special case of
the Maximum Covering Location Problem was defined in response to a
bio-terrorism attack \citep{murali2012facility}. Other examples are
the family of the Obnoxious Facility Location Problems
\citep{church2022review,salazar2025efficient}, where the facilities
to locate have a negative impact on the environment or in the
communities around them. Therefore, the aim in this problem is to
avoid the location of these facilities close to populations. An
example of this family is the well-studied Multiple Obnoxious
Facility Location Problem (MOFLP)
\citep{drezner2019planar,salazar2024basic}.

Usually, problems of these families are formulated on a graph, where
nodes represent demand points, candidate facility sites, or both,
and edges capture the underlying network structure such as
transportation routes or communication links. This graph-based
representation allows distances, costs, and capacities to be naturally
modelled, providing a mathematical framework for analysing and
solving location-allocation decisions. Using graph theory, facility
location models can capture the complexity of real-world systems,
from logistics and supply chains to telecommunication and public
service networks.

Recently, many facility location problems have considered relaxing
some edges of the graph or reducing the demand of some nodes to
obtain better solutions. This relaxation is known as upgrading.
Some examples of upgrading problems are the \textit{p}-Center Problem
with Upgrades \citep{anton2023discrete} and the Maximal Covering
Problem with Upgrades \citep{baldomero2024complexity}. In both cases,
classical FLPs have been modified to allow the possibility to vary
some characteristics of the graph.

In this paper, the Induced \textit{p}-Median Problem with Upgrading
(IpMU) \citep{espejo2023p} is studied. As shown in the previous work,
IpMU is based on the optimization of e-commerce logistics. Beyond
this, IpMU can also represent critical decisions in urban transportation
and public service networks. For instance, in public healthcare
systems, authorities must decide where to locate emergency medical
centers while investing in road or communication infrastructure to
reduce response times. Upgrading certain road segments or improving
digital coordination between facilities can significantly decrease
effective transportation costs, leading to faster patient transfers
and better use of limited budgets.

Formally, the IpMU consists of obtaining the location of $p$ medians
in a bi-network graph and distributing a budget $B$ intended to
upgrade some arcs of the network. The most remarkable feature of
this problem is the bi-network characteristic, which is associated
to the real-world scenario where the time of travel and the cost
of travel are explicitly separated magnitudes. This separation is
defined with two different values in each arc of the input graph.
The set of selected medians will serve the rest of the customers
in the network, always following the fastest path according to the
time of travel. However, the objective function aims to minimize
the cost of travel, considering that upgrading can be applied on
edges to reduce travel cost. As far as we know, the IpMU is the
first study to address $p$-median location problem that combines
arc upgrading with dual-weight networks.

The state of the art of IpMU is so far completely based on mathematical
programming models \citep{espejo2023p}. However, in this paper, we
propose a novel approach based on metaheuristic methods. Metaheuristic
algorithms emerged in the 1980s as an alternative to classical
mathematical programming algorithms to solve optimization problems
due to the long execution times that they require when the size of
inputs grows \citep{sorensen2013metaheuristics}. These methods have
proven their efficiency in different areas of operation research
such as clustering problems \citep{martin2022strategic}, arrangement
problems \citep{robles2025multi}, and also FLP \citep{herran2020parallel}.

State-of-the-art methods based on the Branch and Bound algorithm
need to calculate variables related to node selection and variables
related to edge-upgrading in the same process. This implies that
the base model obtained has undesirable characteristics for its
resolution, such as non-linearity or a high number of variables,
that cause large execution times. In this work, a metaheuristic
algorithm is proposed in a two-phase resolution approach, attaching
the benefits of heuristic and exact methods. Firstly, the selection
of the medians is performed using a method based on the well-known
Greedy Randomized Adaptive Search Procedure (GRASP) metaheuristic
\citep{feo1995greedy,resenderibeiro2016GRASPbook}. The selection
of the medians without attending the edge-upgrading sub-problem
allows the algorithm to reduce the computational time and also,
enables the separation of the subrogated model for the subproblem.

In the study of the subrogated model for the edge-upgrading subproblem, we propose a linear formulation with fewer variables than the original one that could be exactly solved in a reduced execution time with commercial solvers. However, a detailed study of this new developed model shows that it could be solved without the necessity of this software, unlike mathematical model-based algorithms, obtaining a polynomial complexity method that accelerates the resolution of the problem, as it will be shown in the experimental results.

We tested our new proposal with the state-of-the-art benchmark composed of a total of 960 instances of different characteristics. However, due to the small size of these instances (a maximum of 80 nodes in the graph), we have extended the benchmark with a new set of 270 instances of size of 100, 200 and 500 nodes. 

The results of the comparison with the state-of-the-art demonstrate the strong performance of the GRASP proposal, achieving 958 optimal solutions out of 960 in the benchmark set spending less than three seconds on average, while the execution time of the state-of-the-art method increases exponentially. In the case of the new benchmark, the results indicate a further improvement in performance, attaining 268 best-known values compared to 148 achieved by the previous proposal, while requiring an execution time approximately two orders of magnitude shorter.

In the experimentation phase, we observed a high variability in execution times for different instances. Therefore, we assessed the complexity of the instances using the trajectories followed by the algorithms. As a result, we include the findings of this process in this paper, which allows us to identify a possible explanation for their different observed behaviour. This process can be used in future comparisons of algorithms using the proposed extended benchmark.

The remainder of the paper is organized as follows. Section ~\ref{sec:sota} includes a literature review of resolution methods of FLP. Section~\ref{sec:problem} provides a formal description of the IpMU. Section~\ref{sec:algorithmic-proposal} presents the algorithmic approach proposed in this work. Section~\ref{sec:Experiments} details the computational experiments, including parameter tuning, comparisons with state-of-the-art methods, and additional analyses of instances. Finally, Section~\ref{sec:conclusions} summarizes the main findings and provides directions for future research.

\section{Literature review}
\label{sec:sota}

The FLP family is extensive and its taxonomy is based on a variety of different characteristics. Taking into account the mathematical characteristics of the problem model, unilevel and multilevel problems (MFLPs) can be found in the literature \citep{ortiz2018multi}, where several objective functions correspond to different constraints in the complete model. Considering the mathematical characteristics, problems are classified as discrete or continuous, depending on whether the domain of solutions, and therefore the search space, is discrete or continuous. In the case of combinatorial optimization, this family includes classical problems such as the \textit{p}-Median problem \citep{alp2003efficient}. In continuous space, other problems can be found, such as the Fermat-Weber problem \citep{cooper1968extension}. In this section we will focus in unilevel combinatorial FLP problems, which is the category of problems the IpMU belongs to.

Typically, algorithms created to solve problems from the FLP family are divided into two well differentiated groups: those based on mathematical resolution, such as Branch and Cut or Branch and Price algorithms, among others; and heuristic and metaheuristic methods, which spend reduced execution times, but they do not always guarantee the optimality of the obtained solutions.

One of the most studied problems in the FLP family is the Uncapacitated Facility Location Problem (UFLP) where several facilities must be opened in order to bring service to the group of communities around them, without considering their capacity. In this case, the classical approach is based on mathematical models, with a primal-dual algorithm \citep{galvao1989method}. However, this kind of processes are normally constrained to small instances due to the larger execution times. The most recent approaches, which conform the state of the art, are based on heuristic evolutionary algorithms, which bring very good results for larger instances \citep{zhang2023fast,sonucc2023adaptive}.

If demand of the communities and capacity of the facilities are considered, we can find another of the most studied FLP problems: the Capacitated Facility Location Problem (CFLP). Again, the classical method to solve the CFLP is based on an exact algorithm, in this case, a Branch and Bound method with Lagrangian relaxation \citep{nauss1978improved}.
In order to deal with more complex instances, metaheuristic proposals emerge, such as GRASP and Path Relinking algorithms conforming the state-of-the-art of the CFLP \citep{albareda2017heuristic,albareda2023some}.

% Changing the point of view we can find the Maximal Covering Location Problem (MCLP). In this case, instead of installing all the facilities necessary to meet total demand, demand should be maximized with a limited number of facilities. The general process of research is also replicated in this case. The literature brings a first phase of research based on mathematical models \citep{downs1996exact} and finally the complexity of the problems shows the efficiency of heuristic methods \citep{revelle2008solving}.

As evolution of these classical problems, authors in the literature have studied the inclusion of different metrics on the same graph, leading to bi-network location problems. On this topic, we can find Cost-Distance Location Problem \citep{meyerson2008cost}, where authors study an heuristic approach for different problems (also facility location) when cost and distance metrics are involved in the network. In addition, we can find the Median-Path Problem \citep{avella2005branch}. This problem involves finding a path through a network that minimizes both metrics considered, accessibility and total path cost.

As mentioned in Section \ref{sec:intro} the upgrading version of the classical problems is one of the relevant topics in FLP research nowadays. One of the first works identified in the literature is the Upgrading \textit{p}-Median Problem on a Path \citep{sepasian2015upgrading}, where the aim of the problem is to modify the weights of the vertices of a graph with a path structure to obtain a new set of \textit{p} medians that minimizes the sum of the weighted distances to all vertices. Then, more complex versions have emerged in literature. For example, the upgrading version of the Maximal Covering Location Problem (Up-MCLP) \citep{baldomero2024complexity}. In this work, authors state the problem and design optimal methods for special types of graphs based on the mathematical properties of trees or path graphs. Later, the same authors proposed a MILP formulation for the resolution of the problem and a metaheuristic approach based on a local search strategy \citep{baldomero2025edge}.

Concerning the problem under study in this work, the IpMU, fits on both of the problem families previously mentioned: bi-network problems and upgrading problems. The state-of-the-art methods por this problem are completely based on mathematical models \citep{espejo2023p}. In particular, a first formal statement of the problem is defined, then this formulation is linearized, and finally a Branch and Bound algorithm is proposed based on the analysis of the formulation. The experimentation process of this work obtains a specific set of constraints that brings the best Branch and Bound process out of the complete set of constraints defined in the problem. However, no heuristic methods have been developed for this problem in the literature, making our work a novel proposal for the IpMU. 

In brief, Table \ref{tab:sota-review} summarizes this review of the different approaches used to solve FLP problems.

\begin{table}[ht]
\centering
\caption{State of the art of FLP problems.}
\label{tab:sota-review}
\resizebox{\textwidth}{!}{

\begin{tabular}{lll}
\toprule
\textbf{Problem} & \textbf{Exact Method} & \textbf{Heuristic/Metaheuristic Method} \\
\midrule \midrule
\multirow{2}{*}{\textit{UFLP}} &
Primal-Dual Method (1989) \citep{galvao1989method} &
Evolutionary Algorithm (2023) \citep{zhang2023fast} \\
& & Parallel Evolutionary Algorithm (2023) \citep{sonucc2023adaptive} \\
\midrule
\multirow{2}{*}{\textit{CFLP}} &
Branch and Bound (1978) \citep{nauss1978improved} &
GRASP (2017) \citep{albareda2017heuristic} \\
& & Path Relinking (2023) \citep{albareda2023some} \\
\midrule
\textit{MCLP} &
Branch and Bound (1996) \citep{downs1996exact} &
Heuristic Concentration (2008) \citep{revelle2008solving} \\
\midrule
\textit{Median-Path} &
Branch and Cut (2005) \citep{avella2005branch} & - \\
\midrule
\textit{UpPMP on Path} &
Exact approach (2015) \citep{sepasian2015upgrading} & - \\
\midrule
\multirow{2}{*}{\textit{\textit{Up-MCLP}}} &
Models in special graphs (2024) \citep{baldomero2024complexity} & Local Search (2025) \citep{baldomero2025edge}\\
&
MILP (2025) \citep{baldomero2025edge} & \\
\midrule
\textit{IpMU} &
Branch and Bound (2023) \citep{espejo2023p} & - \\
\bottomrule
\end{tabular}

}
\end{table}
As seen in the literature, neither heuristic nor metaheuristic approaches were proposed for the IpMU problem, as in other problems from the FLP family. Therefore, we aim to use metaheuristics to deal with large real-world scenarios where the mathematical approaches require huge amounts of time or are unable to generate any solution.

\section{Formal description}
% Revisada
\label{sec:problem}
In this section, a formal description of the IpMU is developed based on the latest state-of-the-art work \citep{espejo2023p}. Let $G = (V,A)$ be a directed graph. Consider the set of clients as the set of nodes in the graph $V=\{1,...,n\}$. Each client $i \in V$ has an associated demand $w_i \geq 0$. For each edge $a \in A$ consider two different values $c^1_a \geq 0$ and $c^2_a \geq 0$. The first parameter, $c^1$, quantifies the time required to traverse this edge. The second parameter, $c^2$, is the cost per unit of transport on this edge.  

The IpMU problem consists in obtaining the fixed-size set of nodes that minimizes the cost of transport to all clients using a shortest path to reach each. Following the original formulation and notation, $V$ is defined as the set of candidate nodes, and $1\leq p < n$ of them will be selected as the set of medians of solution $S=\{s_1,...,s_p\}$.

Once $S$ is established, every user has to be served by the median that can reach it in the shortest time. Formally, let $\mathit{FP}(i,j) \subset A$ be the shortest path from $j$ to $i$ according to the $c^1$ metric. The aggregated time of $\mathit{FP}(i,j)$ is defined as $$C^1_{i,j} = \sum \limits_{a\in \mathit{FP}(i,j)} c^1_a.$$ Then, every client $i \in V$ will be served by the median $s^\star \in S$ such that $$C^1_{s^\star,i} \leq C^1_{s,i}, \forall s\in S.$$ Following the same structure, the cost of transport is defined by $$C^2_{i,j} = \sum\limits_{a\in \mathit{FP}(i,j)}c^2_a.$$ If a client $i\in V$ can be served by several medians, the one with the lowest associated value $C^2_{s^\star,i}$ will be chosen.

Finally, the value $B \geq 0$ is defined as the budget to improve the cost of transport of the edge. That means it can be used to reduce the value of $c^2$ of some edges of the graph. The relaxation of these values is limited on every edge by the amount $u_a$.

A solution for the IpMU consists in identifying $p$ medians and distributing the budget $B$ among the arcs of the graph in such a way that it minimizes the sum of the transport costs from the medians to their respective customers, weighted by the requested demand.

Equation \ref{eq:obj} formally formulates the problem as a mixed non-linear optimisation problem. For this proposal, it is necessary to define the following variables:
\begin{itemize}
{
    \item $
    y_j = \begin{cases}
    1, \quad &\text{if node $j$ is established as median.}\\
    0, \quad &\text{otherwise.}
    \end{cases}
    $
    \item $
    x_{ij} = \begin{cases}
    1, \quad &\text{if node $j$ is the closest median to client $i$.}\\
    0, \quad &\text{otherwise.}
    \end{cases}
    $
    }
    \item $b_a \in [0,u_a]$ : The cost reduction applied to edge $a\in A$.
\end{itemize}

\begin{align}
\min \quad & \sum_{i\in V} \omega_i \sum_{j\in V} \left(C^2_{i,j} - \sum_{a\in \mathit{FP}(i,j)} b_a \right) x_{ij} \label{eq:obj} \\
\text{s.t.} \quad 
& x_{ij} \leq y_j, && \forall i \neq j \in V, \tag{1.1} \label{eq:objR1}\\
& \sum_{j\in V} x_{ij} = 1, && \forall i \in V, \tag{1.2} \label{eq:objR2}\\
& \sum_{j\in V} y_j = p, && \tag{1.3} \label{eq:objR3}\\
& y_j + \sum_{\substack{s\in V:\\ C^1_{is} > C^1_{ij}}} x_{is} \leq 1, && \forall i,j\in V, \tag{1.4} \label{eq:objR4}\\
& \sum_{a\in A} b_a \leq B, && \tag{1.5} \label{eq:objR5}\\
& b_a \leq u_a, && \forall a \in A, \tag{1.6} \label{eq:objR6}\\
& b_a \geq 0, && \forall a \in A, \tag{1.7} \label{eq:objR7}\\
& x_{ij} \in \{0,1\}, && \forall i,j \in V. \tag{1.8} \label{eq:objR8} \\
& {y_{j} \in \{0,1\}}, && \forall i,j \in V. \tag{1.9} \label{eq:objR9}
\end{align}

\begin{figure}[ht]
    \centering
    \includegraphics[width=0.60\linewidth]{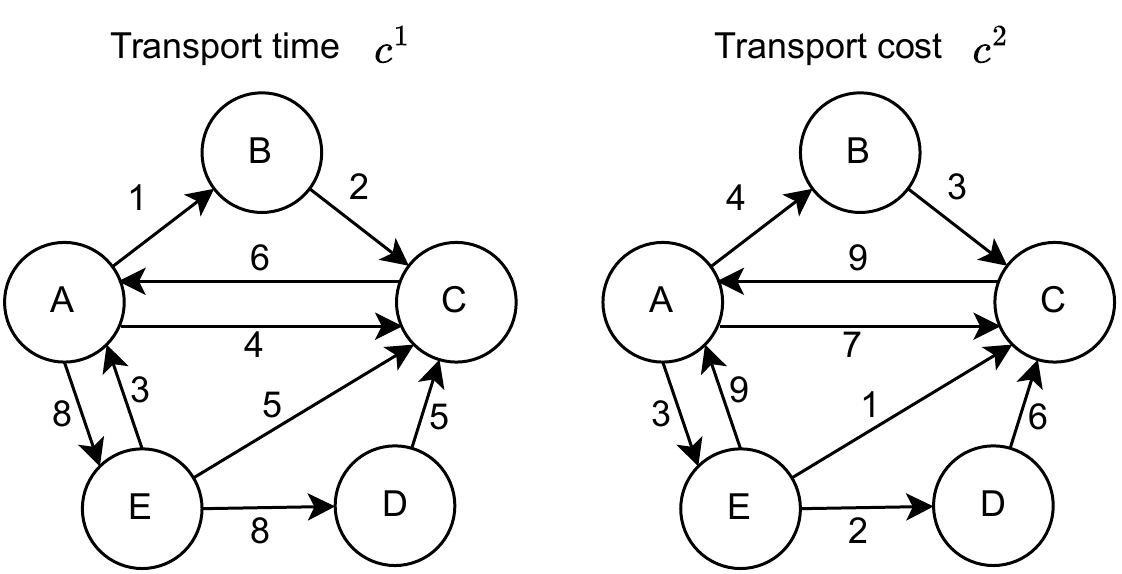}
    \caption{Example instance for the IpMU. It has a budget of $B=2$, a capacity to improve the edges of $u_a=c^2_a$, a uniform demand of $\omega_i=1$ and $p=2$ medians have to be established.}
    \label{fig:ejemplo}
\end{figure}

Figure \ref{fig:ejemplo} shows an example instance of the IpMU, where $c^1$ and $c^2$ values are represented in different pictures of the graph to favour the visualisation. In this case, a budget $B=2$ is considered, and a maximum relaxation value $u_a = c^2_a$ for every edge of the graph. Also, a uniform demand of $\omega_i=1$ for every client is considered. The number of medians to select is $p=2$. Finally, although the nodes of the graph should be numbered from 1 to n, we choose to name them with letters to clarify the example.

\begin{figure}
    \centering
    \includegraphics[width=0.60\linewidth]{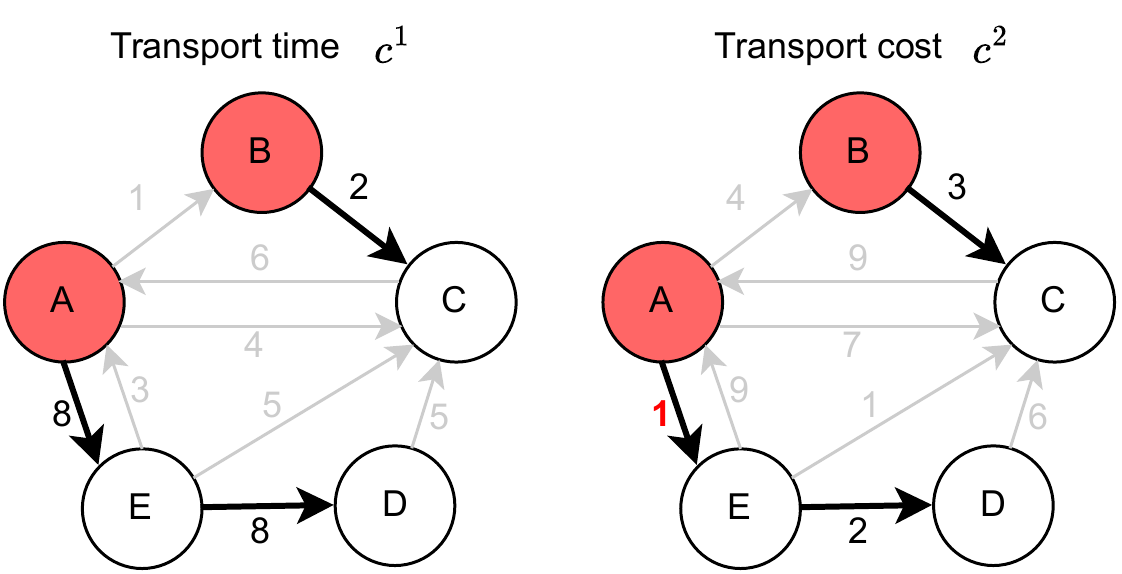}
    \caption{Solution for example instance depicted in Figure \ref{fig:ejemplo}}
    \label{fig:ejemploSolved}
\end{figure}

Figure \ref{fig:ejemploSolved} shows the optimal solution for the example depicted in Figure \ref{fig:ejemplo}. As can be seen, nodes $A$ and $B$ are selected as the medians of the solution. In this case, the fastest way to reach node $C$ is following the edge through $B$, while the fastest paths to $D$ and $E$ start in $A$. See that $\mathit{FP}( {D,A})$ and $\mathit{FP}( {E,A})$ share edge $(A,E)$. This feature makes edge $(A,E)$ the best candidate to be upgraded, relaxing the $c^2$ value from $3$ to $1$. Calculating the cost of these paths, the objective function value obtained is $3+1+3 = 7$ units.

\section{Algorithm proposal}
\label{sec:algorithmic-proposal}
State-of-the-art methods propose linearizations of the model described in Equation \eqref{eq:obj}. In the most recent work \citep{espejo2023p}, three different linear models with different constraints are proposed that improve their performance. Among them, the model named \textit{FL1} is the best proposal  {of the work presented for IpMU}.

The algorithm proposal in our work is based on metaheuristic
procedures \citep{abdel2018metaheuristic}. 
In this case, the proposed
algorithm is based on the popular Greedy Randomized Adaptive Search
Procedure (GRASP) \citep{feo1995greedy,Perez-Pelo2023}. 
In this work,
we proposed a two-phase algorithm to obtain solutions, where the
problem of finding the medians and the problem of upgrading the
edges are detached. This allows us to obtain a reduced and purely
combinatorial search space in the median selection problem, avoiding
the continuous domains of $b_a$, obtaining a faster approach than
a complete resolution.

In this section, we will first define the edge-upgrading problem,
which is obtained by a fast procedure that computes the optimal
relaxation of the costs of the graph for a given set of medians.
Then, the proposal for the median problem based on GRASP will be
described.

\subsection{Solving the edge-upgrading subproblem.}
\label{sub:upgradingsubproblem}

The edge-upgrading sub-problem is described next. Its solution will be incorporated into the complete proposal for the IpMU problem. In different steps of the algorithm, a partial or complete solution must be evaluated. For this purpose, the relaxation of the edges has to be established for a partial or a complete set of medians. Establishing the precise values of each $b_a$ is the \textit{edge-upgrading subproblem}.

Let $S' = \{s_1,...,s_q\}$ be the set of medians selected such that $q \leq p$. In case a partial solution must be evaluated ($q<p$), its quality is given by the objective function of the problem, restricted to select $q$ medians.

Once a set of medians is determined, following the model of Equation \eqref{eq:obj}, the variables $x_{ij}$ and $y_j$ are established, while the variables $b_a$ ($a \in A$) are not.   Under these conditions, the values $\mathit{FP}(i,j), C_{ij}^1$ and $C_{ij}^2$ are constant, because they depend only on $S'$. Following this idea, we can define the IpMU problem subordinated to $S'$ in Equation \eqref{eq:sub}. The definition of this model uses the values mentioned before and could be precalculated when the medians are established. Notice that the model described in Equation \eqref{eq:sub} does not depend on the value of $p$. Then, the model can be used to obtain the optimal relaxation for a complete or a partial solution.

\begin{align}
\max \quad 
& \sum_{i\in V} \sum_{\substack{j\in V \\ x_{ij}=1}} \sum_{a\in \mathit{FP}(i,j)} \omega_i b_a 
\label{eq:sub} \\
\text{s.t.} \quad 
& \sum_{a\in A} b_a \leq B, && \tag{2.1} \label{eq:subR1} \\
& b_a \leq u_a, && \forall a \in A, \tag{2.2} \label{eq:subR2} \\
& b_a \geq 0, && \forall a \in A. \tag{2.3} \label{eq:subR3}
\end{align}

As seen, following this idea, the original formulation is transformed into a linear model with $|A|$ variables and $2|A| + 1$ constraints that result in a simplification considering that the original model in Equation \eqref{eq:obj} minimizes the cost, while the model in Equation \eqref{eq:sub} maximizes the improvement due to the application of a particular budget. As can be seen in this model, each $b_a$ variable is weighted by the sum of the $\omega_i$ values corresponding to the edges $a$ on the shortest path from the median $j$ to customer $i$. Then, we can obtain the corresponding nodes for each of those edges and calculate the sum of their demands. 

Formally, let $N(a)$, described in Equation \eqref{eq:Na}, be the set of nodes of the graph such that $a$ is part of some shortest path.

\begin{equation} 
\label{eq:Na}
\begin{array}{c}
N(a) = \{i \in V \;:\; \exists j \in S\textcolor{red}{'}\;\; s.t. \;( a\in \mathit{FP}(i,j) )\land (C_{ij}^1 = \min\{C_{ik}^1 \;|\; k\in S\textcolor{red}{'}\}) \}
\end{array}
\end{equation}

In case that two different medians $j_1$ and $j_2$ have the same shortest path distance from customer $i$, only the path with the lowest $C^2$ value will be considered. Then we can define the weighted value as in Equation \eqref{eq:wa},

\begin{equation} 
\label{eq:wa}
\begin{array}{c}
\widetilde{\omega}(a) = \sum\limits_{i \in N(a)} \omega_i .
\end{array}
\end{equation}
Finally, the model described in Equation \eqref{eq:sub} is analogous to the model described in Equation \eqref{eq:subRe},

\begin{align}
\max \quad 
& \sum_{a\in A} \widetilde{\omega}(a) \, b_a 
\label{eq:subRe} \\
\text{s.t.} \quad 
& \sum_{a\in A} b_a \leq B, && \tag{5.1} \label{eq:subReR1} \\
& b_a \leq u_a, && \forall a \in A, \tag{5.2} \label{eq:subReR2} \\
& b_a \geq 0, && \forall a \in A. \tag{5.3} \label{eq:subReR3}
\end{align}

It can be seen that the edge-upgrading subproblem is formulated as a Fractional Knapsack Problem \citep{pisinger1998knapsack}, with added bounds, and it is well known that this problem can be efficiently solved by a greedy algorithm. Therefore, we propose a greedy approach to solve this subproblem.

Algorithm \ref{alg:relaxEdges} shows the pseudocode for solving the edge-upgrading subproblem. We name this method \texttt{relaxEdges}. As seen in the pseudocode, the remaining budget $B'$ is first set to the total budget $B$ at step \ref{alg:relaxEdges:b0} and the relaxation values of each edge are set to $0$ in step \ref{alg:relaxEdges:startB}. Then, all the edges $A$ of the graph of the instance are stored and sorted from highest to lowest in $\widehat{A}$ following the weights $\widetilde{\omega}$ received as a parameter and corresponding to the benefits per unit obtained with the relaxation. Then the relaxation value of each edge is established, following the order established in $\widehat{A}$, by making the best use of the remaining budget without exceeding the established limit by each edge $u_a$ (steps \ref{alg:relaxEdges:for} to \ref{alg:relaxEdges:endfor}). Finally, the upgrading values $b$ are returned.

\begin{algorithm}[!htp]
\begin{algorithmic}[1]
    \STATE $B' \gets B$ \label{alg:relaxEdges:b0}
    \STATE $b \gets \{0\}_{|A|}$ \label{alg:relaxEdges:startB}
    \STATE{$\widehat{A}
 \gets \mathtt{sortByWeight} (A,\widetilde{\omega})$} \label{alg:relaxEdges:sortEdges}
    \FOR{$a \in \widehat{A}
$} \label{alg:relaxEdges:for}
        \STATE $b_a \gets \min(B',u_a)$
        \STATE $B' \gets B' - b_a$
    \ENDFOR \label{alg:relaxEdges:endfor}
    \RETURN $b $
\end{algorithmic}
\caption{\texttt{relaxEdges}$(\widetilde{\omega})$}
\label{alg:relaxEdges}
\end{algorithm}
Note that the complexity of Algorithm \ref{alg:relaxEdges} is $\mathcal{O}(|A| \cdot \log |A|)$, derived from the sorting procedure in step \ref{alg:relaxEdges:sortEdges}. The other operation can be performed with linear or constant-time complexity.

Once an algorithm to solve the edge-upgrading subproblem is obtained, an efficient procedure to obtain the evaluation of a set of medians is easily designed. Algorithm~\ref{alg:evaluate} shows the procedure to evaluate a complete or partial solution $S\textcolor{red}{'}$. We denote this process as $\mathcal{F}(S)$. 

In step \ref{alg:evaluate:s0} of the algorithm, the matrix of assignments $x$  {and vector $y$ are} initialized to 0. Subsequently, for each chosen median $j$, the variable $y_j$ is set to 1 and the variables $x_{ij}$ are established so that the median $j$ is closest to node $i$ (steps \ref{alg:evaluate:xjj} to \ref{alg:evaluate:endfor}). In step \ref{alg:evaluate:omega} the weight values are obtained following Equation \eqref{eq:wa}. Then, the vector of real values $b$ is obtained using Algorithm \ref{alg:relaxEdges} in step \ref{alg:evaluate:b}, where $b_a$ represents the relaxation imposed on the edge $a$ of the cost graph. Finally, once all the variables have been set, the value of the objective function is returned. 

\begin{algorithm}[!htp]
\begin{algorithmic}[1]
    \STATE $x \gets \{0\}_{n\times n}$ 
    \label{alg:evaluate:s0}
    \STATE  {$y \gets \{0\}_{n}$} 
    %\STATE $x_{ij} \gets 0 \; \forall i,j \in V$
    \STATE $y_j \gets 1 \; \forall j \in S$ \label{alg:evaluate:xjj}
    \FOR{$j \in S$} \label{alg:evaluate:for}
        \STATE $x_{ij} \gets 1\; \forall i \in \{i\in V : C^1_{i,j} \leq C^1_{i,k}\;  \forall k \in S\}$
    \ENDFOR \label{alg:evaluate:endfor}
    \STATE $\widetilde{\omega} \gets $\texttt{calculateWeights(}$S$\texttt{)} \label{alg:evaluate:omega}
    \STATE $b \gets $ \texttt{relaxEdges(}$\widetilde{\omega}$\texttt{)} \label{alg:evaluate:b}
    \RETURN $\sum\limits_{i\in V} \omega_i \sum\limits_{j\in V} (C^2_{i,j} - \sum\limits_{a\in \mathit{FP}(i,j)} b_a)x_{ij} $
\end{algorithmic}
\caption{$\mathcal{F}(S)$}
\label{alg:evaluate}
\end{algorithm}

Notice that combinatorial search is drastically reduced since the edge-upgrading subproblem is solved with the proposed greedy approach. The metaheuristic process is therefore only limited to selecting $p$ medians from a set of $n$ nodes, and does not need to explicitly consider the $b_a$ variables. For this reason, GRASP procedure only explores a search space of $\binom{n}{p}$ solutions. Although this space grows exponentially, it is significantly smaller than that of the complete IpMU problem and corresponds to a purely combinatorial optimization problem, since the edge-upgrading subproblem is treated separately as a continuous optimization problem.

\subsection{Kuehn-Hamburger heuristic algorithm}
\label{sub:KHAlgo}

Kuehn-Hamburger heuristic \citep{kuehn1963heuristic}, also known as myopic or greedy heuristic, is used as baseline heuristic method in many FLP papers from the literature \citep{ResendeWerneck2004,reese2006solution,mladenovic2007p}.

This method builds a solution incrementally, starting from an empty set of facilities and adding one facility at a time until the desired number $p$ is reached. In each step, the algorithm selects the best candidate to locate the facility between all possible open nodes of the graph. The method locally optimizes the value at each step, without attending to the previous decisions taken at other moments of the process.

Algorithm \ref{alg:KH} shows the pseudocode of the Kuehn-Hamburger heuristic algorithm. See that the implementation of this method in the IpMU problem is not naive because it needs to include our proposal in order to solve the edge-upgrading subproblem.

\begin{algorithm}[!htp]
\begin{algorithmic}[1]
    \STATE $S \gets \emptyset$ \label{alg:KH:S0}
    \WHILE{$|S| < p$}
        \STATE $s \gets \arg \min\limits_{i \in V} \mathcal{F}(S\cup \{i\})$ \label{alg:KH:fmin}
        \STATE$S \gets S \cup \{s\}$ \label{alg:KH:add}
    \ENDWHILE
    \RETURN $S$ \label{alg:KH:end}
 \end{algorithmic}
 \caption{$\mathtt{KH}()$}
 \label{alg:KH}
\end{algorithm}

Step \ref{alg:KH:S0} of the algorithm initializes the current solution as an empty set. Then, in every iteration of the main loop, the algorithm evaluates the contribution of each candidate node to the objective function (step \ref{alg:KH:fmin}), and adds to the solution the best candidate in step \ref{alg:KH:add}. Finally, when the size of the solution reaches $p$ elements, current solution $S$ is returned in step \ref{alg:KH:end}.

This algorithm will be used as a baseline comparison in Section \ref{sub:comparison}, where the GRASP proposal is also compared with the state of the art.

\subsection{GRASP proposal for median selection problem.}
\label{sub:GRASP}

The proposal for solving the IpMU is a Greedy Randomized Adaptive Search Procedure (GRASP) \citep{feo1995greedy,resenderibeiro2016GRASPbook}. This algorithm consists of a randomised-greedy construction phase followed by an improvement phase that has been implemented through a local search.

\subsubsection{Construction phase}

The construction phase consists of the customary greedy randomized construction strategy of GRASP. In this case, the greedy function selected for this method is the same as the objective function. As seen in Section \ref{sub:upgradingsubproblem} this evaluation can be performed correctly for partial or complete solutions.

The algorithm \ref{alg:build} shows the pseudocode for the algorithm construction phase. In step \ref{alg:build:S0} the solution is initialized empty. In step \ref{alg:build:CL0} the candidate list ($\mathit{CL}$) is initialized with all the nodes of the graph. Next, the main loop starts. While the solution is not complete, we obtain the maximum and minimum values of the objective function, $f_{max}$ and $f_{min}$, respectively, obtained after adding each node from $\mathit{CL}$ to the set of medians (steps \ref{alg:build:fmax} and \ref{alg:build:fmin}). With these, a threshold value $\mu$ is obtained using the parameter $\alpha$ that moderates the greediness of the construction method (step \ref{alg:build:th}). In step \ref{alg:build:RCL}, the restricted candidate list ($\mathit{RCL}$) is obtained as the set of candidate nodes that obtain a lower value than $\mu$ when they are evaluated with the actual partial solution. One of these nodes is randomly selected as a median for the solution (steps \ref{alg:build:rnd} and \ref{alg:build:add}) and removed from the candidate list in step \ref{alg:build:minus}. Once the $p$ medians are selected, the obtained solution is returned.

\begin{algorithm}[!htp]
\begin{algorithmic}[1]
    \STATE $S \gets \emptyset$ \label{alg:build:S0}
    \STATE$CL \gets V$  \label{alg:build:CL0}
    \WHILE{$|S| < p$}
        \STATE $f_{max} \gets \max\limits_{i \in CL} \mathcal{F}(S\cup \{i\})$ \label{alg:build:fmax}
        \STATE $f_{min} \gets \min\limits_{i \in CL} \mathcal{F}(S\cup \{i\})$ \label{alg:build:fmin}
        \STATE$\mu \gets f_{max} + \alpha \cdot (f_{min} - f_{max})$  \label{alg:build:th}
        \STATE$RCL \gets \{i \in CL: \mathcal{F}(S\cup \{i\}) \leq \mu\}$  \label{alg:build:RCL}
        \STATE$s \gets random(RCL)$ \label{alg:build:rnd}
        \STATE$S \gets S \cup \{s\}$ \label{alg:build:add}
        \STATE$CL \gets CL\setminus\{s\}$ \label{alg:build:minus}
    \ENDWHILE
    \RETURN $S$
 \end{algorithmic}
 \caption{$\mathtt{Build}(\alpha)$}
 \label{alg:build}
\end{algorithm}

\subsubsection{Improvement phase}
The GRASP improvement phase is implemented as a local search. This local search is defined by a swap move, consisting of exchanging a selected median with an unselected one. Therefore, the method traverses the space of solutions that differ in a median with the initial solution. 

Formally, we can define the Hamming distance \citep{Hamming1950} between two solutions as in Equation \eqref{eq:Hamm}, measuring the number of different elements in both solutions,
\begin{equation} 
\label{eq:Hamm}
\begin{array}{c}
H(S,S') = \frac{|S \; \Delta \; S'|}{2} = \frac{|(S\setminus S') \cup (S' \setminus S)|}{2}.
\end{array}
\end{equation}

Taking into account the previous definitions, given a solution, the local search identifies a better solution among the set of solutions at a Hamming distance of 1.  It stops the search when no such improving solution is available in the swap neighborhood.

Two different strategies are commonly studied in the literature in the implementation of local search \citep{hansen2006first}. The \textit{best improvement} strategy always selects the best option in the complete space of solutions traversed by the local search. On the other hand, the \textit{first improvement} strategy selects the first solution which improves the current solution. For the sake of space, we detail our proposal of \textit{best improvement} since for our problem, as will be mentioned in Section \ref{sub:paramtun}, it obtains better results than the \textit{first improvement} strategy.

Algorithm \ref{alg:SWAP_LS} shows the pseudocode of the \textit{best improvement} local search. In the beginning, the current best solution $S^\star$ is initialized with the incoming solution $S$. Then, a flag variable \textit{improve} is set to \textit{true} to identify when a local minimum has been reached (step \ref{alg:SWAP_LS:flagstart}). The main loop starts iterating while an improvement is made. In step \ref{alg:SWAP_LS:space} the solutions to be explored by the local search are stored in $\Omega$, obtained as a set of complete solutions that differ by one element from the incumbent solution. Then $S'$ is obtained as the highest-quality solution from $\Omega$ in step \ref{alg:SWAP_LS:solution}. If $S'$ is better than the incumbent, this is updated and the search continues (steps \ref{alg:SWAP_LS:ifml} to \ref{alg:SWAP_LS:endifml}). If a local optimum is reached, the solution does not change and the loop is completed. Finally, the incumbent solution $S^\star$ is returned.

\begin{algorithm}[!htp]
\begin{algorithmic}[1]
\STATE $S^\star \gets S$\label{alg:SWAP_LS:start}
\STATE $improve \gets true$ \label{alg:SWAP_LS:flagstart}
\WHILE{$improve$}
    \STATE $improve \gets false$
    \STATE $\Omega \gets \{S \subset V \; :\; (|S| = p) \;\land\; (H(S,S^\star) = 1)\}$ \label{alg:SWAP_LS:space}
    \STATE $S' \gets \arg \max\limits_{S\in \Omega} \mathcal{F}(S)$ \label{alg:SWAP_LS:solution}
    \IF{$\mathcal{F}(S') < \mathcal{F}(S^\star)$}\label{alg:SWAP_LS:ifml}
        \STATE$improve \gets true$ \label{alg:SWAP_LS:flagchange}
        \STATE $S^\star \gets S'$ \label{alg:SWAP_LS:change}
    \ENDIF\label{alg:SWAP_LS:endifml}
\ENDWHILE
\RETURN $S^\star$
\end{algorithmic}
\caption{$\mathtt{LocalSearch}(S)$}
\label{alg:SWAP_LS}
\end{algorithm}

\subsubsection{Final proposal}
Our final GRASP proposal connects both the construction and improvement phases in a multi-start schema. Often, the termination condition of GRASP is defined by a maximum number of iterations. In this case, due to the complexity of the problem, we define the termination condition in relation to two parameters: maximum number of iterations and maximum number of iterations without improvement. This way, we allow the algorithm to finish prematurely in case a given maximum number of iterations without improving is reached.

Algorithm \ref{alg:GRASP} shows the pseudocode for the final GRASP proposal of this work, which receives the parameter $\alpha$ used in the construction phase (see Algorithm \ref{alg:build}), and $maxIters$ and $maxIters_{wi}$, corresponding to the maximum number of iterations and the maximum number of iterations without improvement. First, initial values are established: $S^\star$ as the best known solution, and $iters$ and $iters_{wi}$, which maintain the current number of iterations of the multi-start schema and the current number of iterations without finding a better solution, respectively (steps \ref{alg:GRASP:start} to \ref{alg:GRASP:start3}). Then, the main loop begins and the build and improve phase will execute while maximum iteration parameters are not reached. In steps \ref{alg:GRASP:itermas} and \ref{alg:GRASP:iterwimas} the actual iterations increase by one. Then, the current solution $S$ is built and improved (steps \ref{alg:GRASP:build} and \ref{alg:GRASP:improve}) following the strategies described in algorithms \ref{alg:build} and \ref{alg:SWAP_LS}, respectively. In case the new solution is better than $S^\star$, the iterations without improvement are restarted and the solution is stored (steps \ref{alg:GRASP:if} to \ref{alg:GRASP:endif}). When the loop ends, the best solution found $S^\star$ is returned.

\begin{algorithm}[!htp]
\begin{algorithmic}[1]
\STATE $S^\star \gets \emptyset$\label{alg:GRASP:start}
\STATE $iters \gets 0$\label{alg:GRASP:start2}
\STATE $iters_{wi} \gets 0$\label{alg:GRASP:start3}
\WHILE{$iters \leq maxIters \;\;\land\;\; iters_{wi} \leq maxIters_{wi}$}
    \STATE $iters \gets iters+1$\label{alg:GRASP:itermas}
    \STATE $iters_{wi} \gets iters_{wi}+1$\label{alg:GRASP:iterwimas}
    \STATE $S \gets \mathtt{Build}(\alpha)$\label{alg:GRASP:build}
    \STATE $S \gets \mathtt{LocalSearch}(S)$\label{alg:GRASP:improve}
    \IF{$\mathcal{F}(S) < \mathcal{F}(S^\star)$} \label{alg:GRASP:if}
        \STATE $iters_{wi} \gets 0$
        \STATE $S^\star \gets S$
    \ENDIF\label{alg:GRASP:endif}
\ENDWHILE
\RETURN $S^\star$
\end{algorithmic}
\caption{$\mathtt{GRASP}(maxIters, maxIters_{wi},\alpha)$}
\label{alg:GRASP}
\end{algorithm}

\section{Experimental Results}
\label{sec:Experiments}

Once the algorithmic approach and the analysis of the problem have been described, we present in this section the thorough experimentation that we have carried out for this problem.

All experiments described here were executed on an AMD EPYC 7282 16-core virtual CPU with 32GB of RAM. To perform a fair comparison, the models described in \citep{espejo2023p} have been reimplemented using the well-known tool \texttt{Gurobi 11.0}. In this case, the \texttt{Gurobi} package of \texttt{Java} has been used to solve the problem running version 11.0. The development of the final algorithm has been driven by the \texttt{MORK 0.20} framework \footnote{\url{https://github.com/mork-optimization/mork}} \citep{MORK2} under the \texttt{Java 21.0.1} programming language.

\subsection{Instances and parameter tuning}
\label{sub:paramtun}
In the algorithm experimentation, all state-of-the-art instances have been considered. These are two distinct sets. The first set, called \textit{P}, considers that there is a correlation between cost and time of transport. This correlation can be checked by evaluating the $R^2$ factor between both values, which is above $0.98$ in all cases. The second set of instances, called \textit{R}, does not consider this correlation, and values have been randomly generated. Both sets share the same characteristics for the remaining metrics: $n\in\{20,40,60,80\}$; $m\in[100,500]$; $p\in\{2,3,4,5\}$ and $B\in\{50,100\}$. It should be mentioned that, for each value of $n$, $m$ takes three different values in the given range. In addition, for each combination of parameters, five different instances were generated. Therefore, these values make up two sets of 480 instances of each type, that is, a total of 960 instances studied from the state-of-the-art. For more details on the characteristics of this set, the reader is directed to the state-of-the-art article \citep{espejo2023p}.

However, these instances do not provide a fair comparison with reality, where networks with a larger number of nodes exist. For this reason, new instances were created, following the same method as in the state-of-the-art.  In this case: $n\in \{100,200, {500}\}$; $m\in \{0.25\gamma,0.5\gamma,0.75\gamma\}$, where $\gamma = n\cdot(n-1)$ that corresponds with the number of edges of a complete directed graph to study instances with low, medium and high density; $p\in\{2,3,4,5,10\}$ and $B=100$. Regarding the weights of the cost graph, the values are selected following the same schema as the state-of-the-art authors: in case of \textit{R} type, values are chosen between 0 to 100 following a uniform distribution; in case of \textit{P} type instances $c^2_a := c^1_a + U$ where $U \sim U(1,1.5) $. For every combination of parameters, three instances were created. These instances were generated for both types, \textit{P} and \textit{R}, obtaining two new sets of 90 instances, therefore adding 270 instances to the complete set of study. 

A first experiment allowed us to determine the best values for the
parameters of the different algorithmic components proposed in this
work. In this experiment, a subset of instances from the complete
set is required to avoid overfitting the setup. In this work, we
obtained a set of representative instances following the procedure
described in \citet{InstSelector}, considering 6 different features
of the instances: number of nodes, number of medians to be selected,
density between nodes and edges, minimum demand of the clients,
maximum demand of the clients, and the budget allowed to the edge
upgrade sub-problem.

In particular, the process requires an exploratory landscape analysis to characterize instances as feature vectors, and then selects subsets of instances by clustering and graph analysis. This ensures a more balanced and diverse instance set, reducing the risk that results are biased toward instances that favour a particular algorithm. We determined to select $15\%$ of the complete set, which results in a number of 172 representative instances.

To obtain the best configuration of the parameter values we used the well-known tool \texttt{irace} \citep{lopez2016irace}. This software is based on the iterated racing procedure, which repeatedly samples candidate configurations, evaluates their performance across problem instances, and statistically eliminates worst configurations through a racing process. By iteratively updating the sampling distributions toward the best-performing configurations, irace efficiently explores the parameter space while avoiding premature convergence.

In this case, we have not included the maximum number of iterations that GRASP can perform (\textit{maxIters}) as a parameter for \texttt{irace}, since naturally it will select a value as large as possible, given that the focus of \texttt{irace} is on the quality of the solutions, rather than the execution time. Therefore, the value of this parameter was manually selected with a preliminary experimental phase in which a convergence analysis was performed by randomly changing the value of $\alpha$ on each iteration to study the balance between these two metrics. As a result, the value of \textit{maxIters} is established to 100.

Once the maximum number of iterations is chosen, a threshold is established and we can delegate the choice of the maximum number of iterations without improvement ($maxIters_{wi}$) to \texttt{irace}, since it will not exceed the \textit{maxIters} value.

Hence, three parameter values are tuned: the value of parameter $\alpha \in [0.00,1.00]$, which measures the greediness of the constructive phase; the strategy of the local search, the best improvement (\textit{BI}) or first improvement (\textit{FI}); and the maximum iterations that the algorithm can perform without improving ($maxIters_{wi} \in \{0,...,100\}$).

The tool \texttt{irace} was executed with the set of representative instances in a total of 10000 experiments. Table \ref{tab:irace} shows the parameters obtained. After the whole execution of \texttt{irace} only one elite configuration survived. This configuration states that the value for $\alpha$ is $0.51$, a best improvement search strategy is selected for the local search, and a maximum of 29 iterations without improvement is set.

\begin{table}[ht]
\centering
\label{tab:irace}
\caption{Parameters values obtained by \texttt{irace}.}
\begin{tabular}{lccc}
\toprule
\textbf{Parameter} & $\alpha$  & LS    & $maxIters_{wi}$ \\ \midrule
\textbf{Value}          & 0.51 & \textit{BI}      & 29         \\
\bottomrule
\end{tabular}

\end{table}

The values obtained by \texttt{irace} conform to the final proposal for this work, which will be compared with the state-of-the-art methods in the next section.

\subsection{Comparison with the state of the art}
\label{sub:comparison}
This section presents a comparison of the results with the state of the art for the IpMU, shaped by \citet{espejo2023p}, where three different mathematical models (\textit{FL1,FL2} and \textit{FL3}) are presented to solve the problem with some extra constraints that accelerate their executions. Since \textit{FL1} obtains the best results reinforced with additional constraints, we have implemented this model that includes those constraints ((12), (14) and (15), according to \citet{espejo2023p}) using \texttt{Gurobi}. To perform a fair comparison, we have executed this model on the same machine as the GRASP. For the sake of clarity, we refer to the implementation of the state-of-the-art formulation as \textit{FL1}. In addition, we include another proposal based on the Kuehn-Hamburguer heuristic (see Section \ref{sub:KHAlgo}). The idea of this proposal is to be used as a baseline for the comparison. For the sake of clarity, we refer to this proposal as \textit{KH} in this section.

To compare the results of the algorithms on different sets of instances, we use tables with the same structure. The first column shows the instance set type, either \textit{P} of \textit{R}. The second column shows the name of the algorithm. The third one, \textit{Avg. F(S)}, shows the average value of the objective function among the corresponding instances. The fourth column, \textit{T. (sec)}, shows the average time in seconds that the algorithm needs to execute these instances. The fifth column, \textit{Dev. (\%)}, shows the average deviation from the best known result that the algorithm obtains among all these instances. Finally, the sixth column, \textit{\# Opt.}, shows the number of optimum values obtained by each algorithm. Notice that, despite \textit{FL1} is able to certify the optimal values when its execution ends, when this algorithm is unable to finish its execution before the time limit of one hour, this column reports both the best known results in each instance for every algorithm, \textit{\# B.}, and the number of verified optimum values in brackets. Similarly, the deviation value is not calculated from the optimal value, but from the best known value.

%The next table shows the aggregate results of all the instances studied; for a detailed comparison of each instance, the reader is referred to the excel file of the appendix.

Table \ref{tab:res-small} shows the aggregated results of \textit{GRASP},  {\textit{KH} and} \textit{FL1} in the benchmark instances proposed by the state-of-the-art \citep{espejo2023p}. 
In the case of instances \textit{P}, \textit{KH} the algorithm only reaches a 10\% of the best knowkn values, which implies a 14\% of deviation in terms of quality. However, it requires a very reduced execution time of 0.04 seconds on average. In the case of \textit{GRASP} and \textit{FL1}, both algorithms are able to obtain the best result in the complete set of 480 instances. The only difference between both algorithms lies in the runtime metric, where \textit{GRASP} obtains an average value of $2.59$ seconds, while \textit{FL1} needs $4.11$ seconds on average.  

On the other hand, in the case of set \textit{R}, \textit{GRASP} obtains all the optimal values of the set of instances, for all but two instances of the 480, corresponding to a negligible average deviation of $0.002\%$. However, the improvement in execution time is remarkable in contrast to \textit{FL1}, which takes $42.50$ seconds on average, while \textit{GRASP} takes just $3.01$ seconds, representing $7.1\%$ of the execution time of \textit{FL1}. In the case of \textit{KH}, a behaviour similar to the one observed for the \textit{P} instances is replicated. The algorithm spends a very short execution time, but this comes at the cost of producing the poorest solution quality among all evaluated methods.

\begin{table}[ht]
\centering
\caption{Comparison results on the state-of-the-art instances benchmark ($n\in \{20,40,60,80\}$).}
\label{tab:res-small}
\begin{tabular}{clrrrr}
\toprule
Type& Algorithm     & Avg. F(S)      & T. (sec) & Dev. (\%) & \# Opt. \\\midrule
\multirow{3}{*}{\textit{P}}& \textit{GRASP} & \textbf{41457.45} & 2.59  & \textbf{0.000\%} & \textbf{480}      \\
&  {\textit{KH}}   &  {47819.62} &  {\textbf{0.04}}  &  {14.021\%} &  {49}\\
& \textit{FL1}   & \textbf{41457.45} & 4.11  & \textbf{0.000\%} & \textbf{480}\\   \midrule  
\multirow{3}{*}{\textit{R}} &\textit{GRASP} & 69278.86 & 3.01  & 0.002\% & 478      \\
& {\textit{KH}}   &  {80851.56} &  {\textbf{0.03}}  &  {18.748\%} &  {43}\\
&\textit{FL1}   & \textbf{69276.87} & 42.50  & \textbf{0.000\%} & \textbf{480}\\   \bottomrule
\end{tabular}

\end{table}

\begin{figure}[ht]
    \centering
    \includegraphics[width=0.85\linewidth]{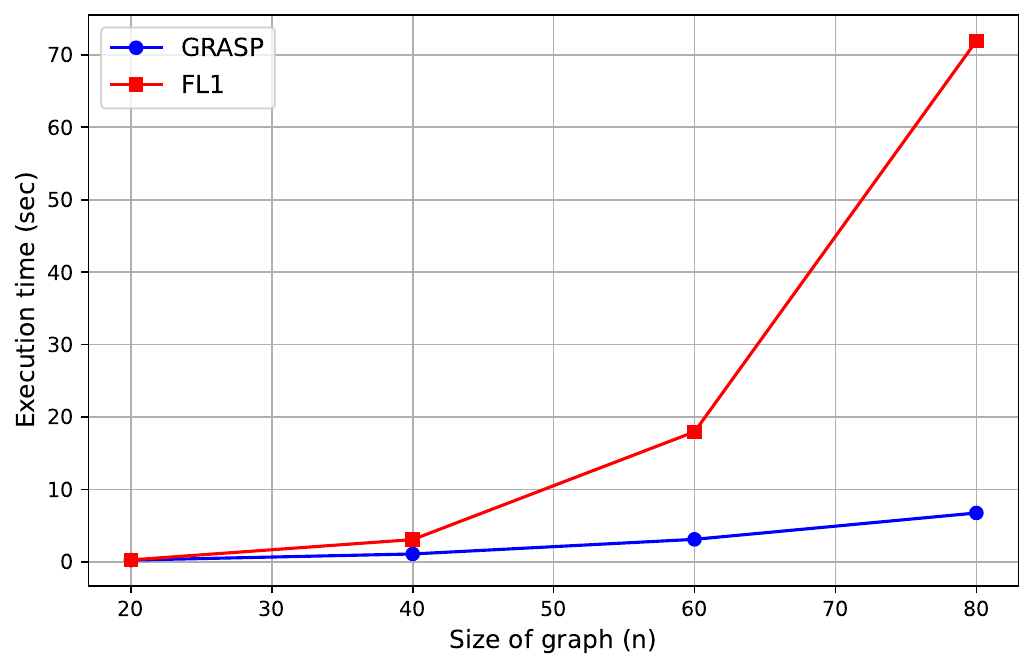}
    \caption{Average execution time comparison between \textit{GRASP} and \textit{FL1} in the set of instances of the state-of-the-art \citep{espejo2023p}.}
    \label{fig:timeComparison}
\end{figure}

A detailed view of the execution time of both methods shows that the performance of the \textit{GRASP} algorithm is notable. Figure \ref{fig:timeComparison} shows the average execution time of each algorithm depending on the size of the instance considering all state-of-the-art instances. As can be seen in the evolution of time, \textit{FL1} follows an exponential-like behaviour that implies long execution times even when the graph size is not too big (80 nodes). However, \textit{GRASP} behaviour is much smoother, maintaining an average execution time shorter than $10$ seconds even in the largest instances studied in the previous work of IpMU. This figure shows that the GRASP algorithm spends smaller execution times in all the graph sizes proposed in the literature benchmark. In smaller size instances GRASP and FL1 times are comparable, however, there is a clear upward trend in execution time for both algorithms, demonstrating the need for another algorithmic proposal as GRASP as the number of nodes in the instances increases. 

As stated before, an additional experiment was conducted using larger instances with sizes of 100 and 200 nodes. It is remarkable that these sizes do not represent an overstatement since there exist networks with a larger number of nodes studied even in \textit{p}-Median problems \citep{garcia2011solving,mu2020solving,ResendeWerneck2004}. In order to maintain reasonable execution time values for these experiments, the execution of \textit{FL1} has been limited to one hour, and the reported value of the objective function is the best value \texttt{Gurobi} obtains in this time.

Tables \ref{tab:res-bigP} and \ref{tab:res-bigR} show the comparison of the proposed instances with sizes $n\in \{100,200\}$. In order to show a detailed view of these results, we append an extra column indicating the size of the analysed instances, denoted as \textit{Size (n)}.

In the case of \textit{P} type instances in Table \ref{tab:res-bigP}, we can observe that \textit{GRASP} and \textit{FL1} obtain all the optimal results for instances with $100$ nodes, but with significant differences in execution time, where our proposal spends just $18.89$ seconds, a third of the time compared to the model. 

In the case of instances with $200$ nodes, the superiority of the \textit{GRASP} proposal is most remarkable, obtaining the best known values for all the 45 instances, and reaching all the 34 certified optimum values. On the other hand, the \textit{FL1} model is not able to obtain the best results in two instances since the time limit of one hour is reached. Furthermore, for these instances of $200$ nodes, the execution time increases significantly. The \textit{FL1} model requires more than 2200 seconds, while \textit{GRASP} needs 121 seconds on average, which implies a reduction in execution time of $94.65\%$.

In the case of \textit{KH} performance, the capacity to achieve good results worsens even further, reaching a deviation of 28\% and 22\% compared to the other proposals. Although the execution time is very short, the results obtained are of such poor quality that the method cannot be considered a viable proposal for the IpMU.

%Finally, notice that our \textit{GRASP} proposal is able to reach all the certified optimum values and, for the 11 instances with size $n=200$ where \texttt{Gurobi} is not able to certify the optimum value, our \textit{GRASP} proposal gets the best results in two more than \textit{FL1}. 

\begin{table}[ht]
\centering
\caption{Comparison results for \textit{P} type instances with $n\in \{100,200\}$.}
\label{tab:res-bigP}
\resizebox{\textwidth}{!}{
\begin{tabular}{cclrrrr}
\toprule
Type & Size (n) & Algorithm & Avg. F(S)      & T. (sec) & Dev. (\%) & \# B. (\# Opt.) \\\midrule
\multirow{6}{*}{\textit{P}}&\multirow{3}{*}{\textit{100}}&\textit{GRASP} & \textbf{24244.56} & 18.89  & \textbf{0.000\%} & \textbf{45 (45)} \\
&&  {\textit{KH}}  &  {30340.02} &  {\textbf{0.08}}  &  {28.110\%} &  {0 (0)}\\
&&\textit{FL1}   & \textbf{24244.56} & 57.67  & \textbf{0.000\%} & \textbf{45 (45)}\\   \cmidrule(r){2-7} 
&\multirow{3}{*}{\textit{200}}&\textit{GRASP} & \textbf{40808.82} & 121.13  & \textbf{0.000\%} & \textbf{45 (34)}  \\
&&  {\textit{KH}}  &  {49240.31} &  {\textbf{0.46}}  &  {22.145\%} &  {0 (0)}\\
&&\textit{FL1}   & 40844.25 & 2264.88  & 0.100\% & 43 (34)\\   \bottomrule
\end{tabular}
}

\end{table}

In the case of \textit{R} type instances, shown in Table \ref{tab:res-bigR}, the performance is comparable. Our \textit{GRASP} proposal reaches all the certified optimum values in the instances with size $n=100$, except one, obtaining a negligible deviation of $0.007\%$. However, these results are obtained spending an average execution time of 24 seconds, while the model requires 492 seconds, which is more than 20 times slower than the metaheuristic proposal.

Analysing the instances with $200$ nodes, the quality difference between both proposals is more remarkable. In this case, the \textit{FL1} model is not able to certify the optimality in any of the instances, obtaining the best known solution in a third of the studied instances. On the other hand, our \textit{GRASP} proposal reaches 44 of the 45 best known values, obtaining a negligible deviation value of $0.014\%$, while the \textit{FL1} model reaches a deviation of $4.467\%$, more than 2 orders of magnitude above our proposal.

In addition to the improvement in quality obtained by the metaheuristic algorithm, there is also a significant improvement in execution time. The average execution time of the \textit{GRASP} proposal represents a saving of $95.5\%$ for instances with $200$ nodes, since \textit{FL1} reaches the time limit of one hour in all of them. Finally, analysing the \textit{KH} results, the values obtained show once again that a more intelligent metaheuristic approach, such as \textit{GRASP}, is required for IpMU.

\begin{table}[ht]
\centering
\caption{Comparison results for \textit{R} type instances with $n\in \{100,200\}$.}
\label{tab:res-bigR}
\resizebox{\textwidth}{!}{
\begin{tabular}{cclrrrr}
\toprule
Type & Size (n) & Algorithm & Avg. F(S)      & T. (sec) & Dev. (\%) & \# B. (\# Opt.) \\\midrule
\multirow{6}{*}{\textit{R}}&\multirow{3}{*}{\textit{100}}&\textit{GRASP} & 379222.09 & 24.01  & 0.007\% & 44 (44)      \\
&&  {\textit{KH}}  &  {502088.25} &  {\textbf{0.08}}  &  {35.627\%} &  {0 (0)}\\
&&\textit{FL1}   & \textbf{379204.78} & 492.50  & \textbf{0.000\%} & \textbf{45(45)}\\   \cmidrule(r){2-7} 
&\multirow{3}{*}{\textit{200}}&\textit{GRASP} & \textbf{943254.88} & 159.39  & \textbf{0.014\%} & \textbf{44} (0)      \\
&&  {\textit{KH}}  &  {1272752.50} &  {\textbf{0.51}}  &  {36.927\%} &  {0 (0)}\\
&&\textit{FL1}   & 981851.37 & 3600.00  & 4.467\% & 15 (0)\\   \bottomrule
\end{tabular}
}
\end{table}

To conclude the benchmark comparison, we analysed instances with size $n=500$. This case is separated from the other proposed instances because the model implemented on Gurobi is not able to solve any of them. The solution process demands a substantial amount of memory, and once the hardware memory limit is exceeded, the Java Virtual Machine (JVM) terminates the process before a feasible solution can be obtained. Table \ref{tab:res500} summarizes the results obtained. As can be seen, GRASP obtains a good solution for all instances in an average time of 30 minutes, without being limited by the available memory.

\begin{table}[ht]
\centering
\caption{Comparison results for instances with $n=500$.}
\label{tab:res500}
\resizebox{\textwidth}{!}{

\begin{tabular}{cclrrrr}
\toprule
Type & Size (n) & Algorithm & Avg. F(S)      & T. (sec) & Dev. (\%) & \# B. (\# Opt.) \\\midrule
\multirow{3}{*}{\textit{P}}&\multirow{3}{*}{\textit{500}}&\textit{GRASP} & 102265.76 & 1576.98  & \textbf{0.000\%} & \textbf{45 (0)}      \\
&& \textit{KH}  & 122482.49 & 7.32  & 20.20\% & 0 (0)\\
&&\textit{FL1}   & - & -  & - & - \\  \midrule 
\multirow{3}{*}{\textit{R}}&\multirow{3}{*}{\textit{500}}&\textit{GRASP} & 3568726.29 & 1539.88  & \textbf{0.000\%} & \textbf{45 (0)}      \\
&& \textit{KH}  & 4602170.51 & 8.43  & 21.65\% & 0 (0)\\
&&\textit{FL1}   & - & -  & - & - \\   \bottomrule
\end{tabular}

}
\end{table}

Therefore, these results demonstrate that the proposed \textit{GRASP} algorithm improves the current state-of-the-art results. For small instances, where mathematical programming approaches usually outperform heuristic methods, \textit{GRASP} achieves solutions of comparable quality in significantly shorter execution times. For larger instances, the differences are even more pronounced in both computation time and objective function value: in the largest case that FL1 obtains solutions with 200 nodes, \textit{FL1} exhibits a 4.47\% deviation while requiring more than 20 times the execution time of the proposed metaheuristic. Notice that the GRASP proposal is not in any case above 0.015\% deviation. Besides, \textit{FL1} is not able to finish its execution in the largest instances due to memory limitations. Again, this shows the superiority of the metaheuristic proposal in terms of execution, where the results obtained are much better using 20 times less time and without problems due to the memory restrictions. 

Finally, to provide a theoretical assessment of the performance
differences among the three algorithms, an additional statistical
analysis was conducted with the bigger instances of the problem
with 100 and 200 nodes since FL1 is not able to solve $n=500$
instances. Specifically, a performance comparison was carried out
using the Bayesian approach proposed in \citet{calvo2019bayesian}.
This method ranks the algorithms according to the probability of
achieving better results than the others and provides corresponding
confidence intervals for these probabilities. The probability of
winning is estimated through a Plackett-Luce model based on the
ranking of the algorithms for each instance. Further details on the
methodology can be found in the referenced paper.

\begin{figure}[ht]
    \centering
    \includegraphics[width=0.8\linewidth]{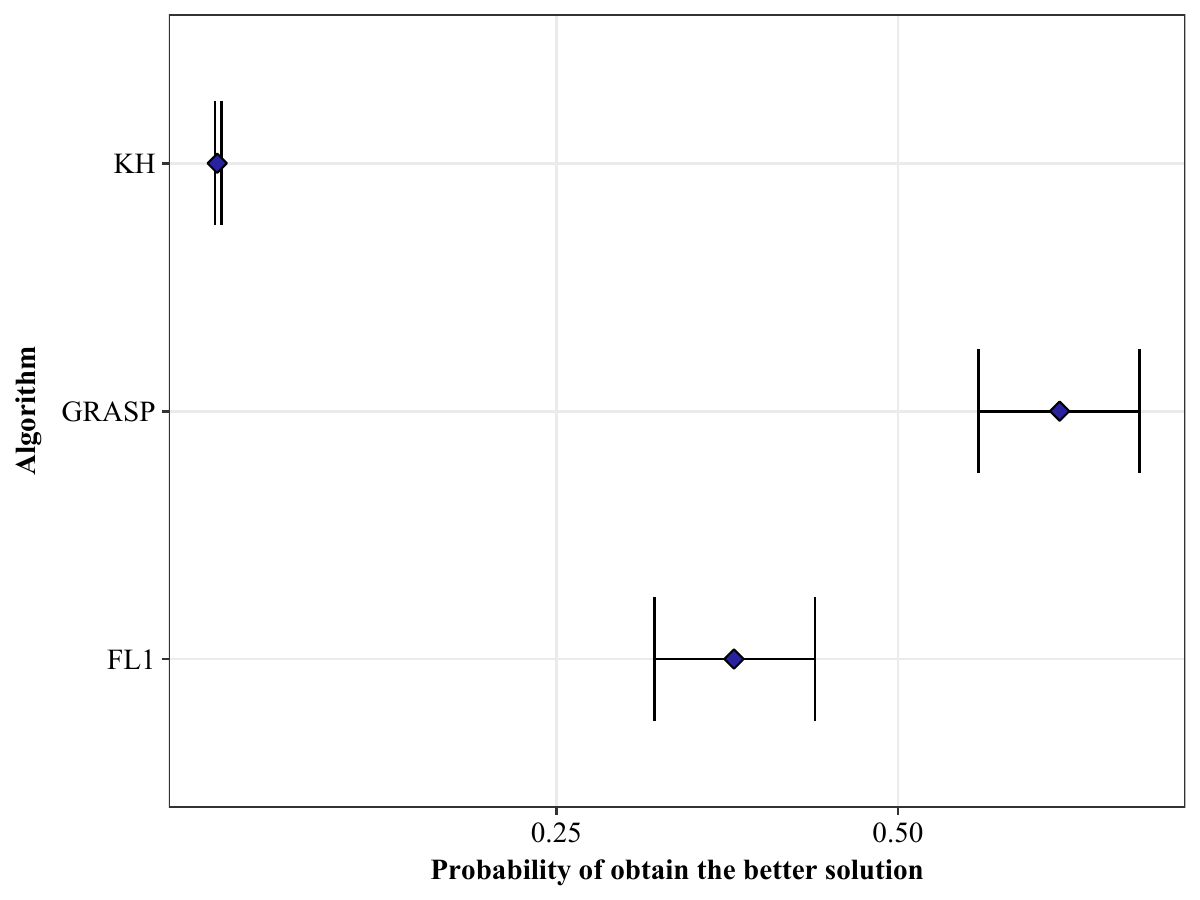}
    \caption{ {Bayesian Inference Statistical Test.}}
    \label{fig:bayesian}
\end{figure}

Figure \ref{fig:bayesian} presents the results of the Bayesian analysis. As shown, the \textit{GRASP} approach achieves a higher probability of outperforming the other methods. Specifically, this probability is $0.618$, whereas the second-best method, \textit{FL1}, attains a value of $0.380$. Moreover, the credible intervals of the algorithms do not overlap with those of the other approaches, indicating that, even in the worst-case scenario, the probability of \textit{GRASP} obtaining a better result remains higher than the current state-of-the-art method. Finally, the KH heuristic obtains the worst performance with a negligible $0.002$ probability of winning.

% \sout{Therefore, these results show that the \textit{GRASP} is able to improve the state of the art results both in solution quality and execution time. In addition, it can be seen that the difference with the \textit{FL1} proposal from the state of the art clearly increases with the size of the instances.}

\subsection{Hardness analysis of the instances}
\label{sub:hard}
As seen in Section \ref{sub:comparison}, solving instances of type \textit{R} requires more execution time for both algorithms, and the deviations over the best known results are also higher, while instances of type \textit{P} do not present such a clear trend. Therefore, it seems plausible that the hardness of the instances could be linked to the correlation between time and cost of transport. 

Typically, visualization of the space of solutions helps in identifying features of the instances that could be related to the performance of the algorithms. This process is usually simple in  some continuous problems. However, combinatorial optimisation problems do not have an easy visualisation of their solution space. 

%This seems to indicate that the search space of the instances with correlation between time and cost of transport is less complex with respect to the instances that this correlation does not exist.

In this paper, we use a method to represent the search space of the problem in order to study the hardness of the instances.

%obtain a confidence measure of how the instances in relation to this feature change, in this article, we have worked on a method to represent the search space of the problem.

In the IpMU problem, a solution is completely represented by the medians selected among the graph nodes because, as mentioned in Section \ref{sub:upgradingsubproblem}, the edge-upgrading subproblem can be optimally solved, without including additional information in the solution. This means that, given an IpMU instance with $n$ nodes and $p$ medians to select, there is a total of $\binom{n}{p}$ different solutions. Although this value grows exponentially, when the parameters are not extremely large, they can be fully enumerated.

As mentioned in Section \ref{sub:GRASP}, the local search proposed for the GRASP method is based on the swap move, that is, between the incoming solution and the next improved solution there is a unique median of difference when the solution changes or no difference when the input solution is a local optimum. In addition, the best improvement strategy implies that this method is totally deterministic. In this case, the local search brings a way to move between solutions of the instance following the trajectory of the search.

Hence, we can define the Search Space Graph ($\mathit{SSG}$) of an instance of the IpMU. The $\mathit{SSG} = (\mathcal{V},\mathcal{E})$ is a directed graph where the set of nodes $\mathcal{V}$ is the set of possible solutions of the instance, that is, the set of subsets of $p$ elements of $V$. On the other hand, the set of edges is formed by the trajectory followed by the local search method. Hence, for every node of $\mathcal{V}$ there exists a unique edge to the solution that minimizes the value of the objective function with a single median change, whose direction follows the search trajectory. These properties are formally stated in equations \eqref{eq:visual-v} and \eqref{eq:visual-v},
\begin{equation}
\label{eq:visual-v}
    \mathcal{V} = \{S \subseteq V : |S| = p\},
\end{equation}
and
\begin{equation}
\label{eq:visual-e}
    \mathcal{E} = \{(S,\arg \min \limits_{\substack{S' \in \mathcal{V} \\ H(S,S')=1 \\ \mathcal{F}(S')<\mathcal{F}(S) } }\mathcal{F}(S')) \; \;\forall S \in \mathcal{V}\} .
\end{equation}
Then, the number of nodes of \textit{SSG} is equal to the number of solutions of the problem $|\mathcal{V}| = \binom{n}{p}$. On the other hand, the number of edges is $|\mathcal{V}|-o^\star$ where $o^\star$ is the number of local optimum solutions of the instance. This value is explained by the fact that for every node only one edge points to the best neighbour that improves its objective function value except if none of them can improve, i.e., is a local optimum. Notice that \textit{SSG} is a graph without cycles, where a local optimum is a node with no outgoing edges. 
This feature implies that \textit{SSG} is a forest graph, understanding a forest as a set of tree graphs, where for every local optimum there exists a disjoint tree with the nodes that converge to a local optimum considering the trajectory of the local search.

\begin{figure}
    \centering
    \includegraphics[width=0.8\linewidth]{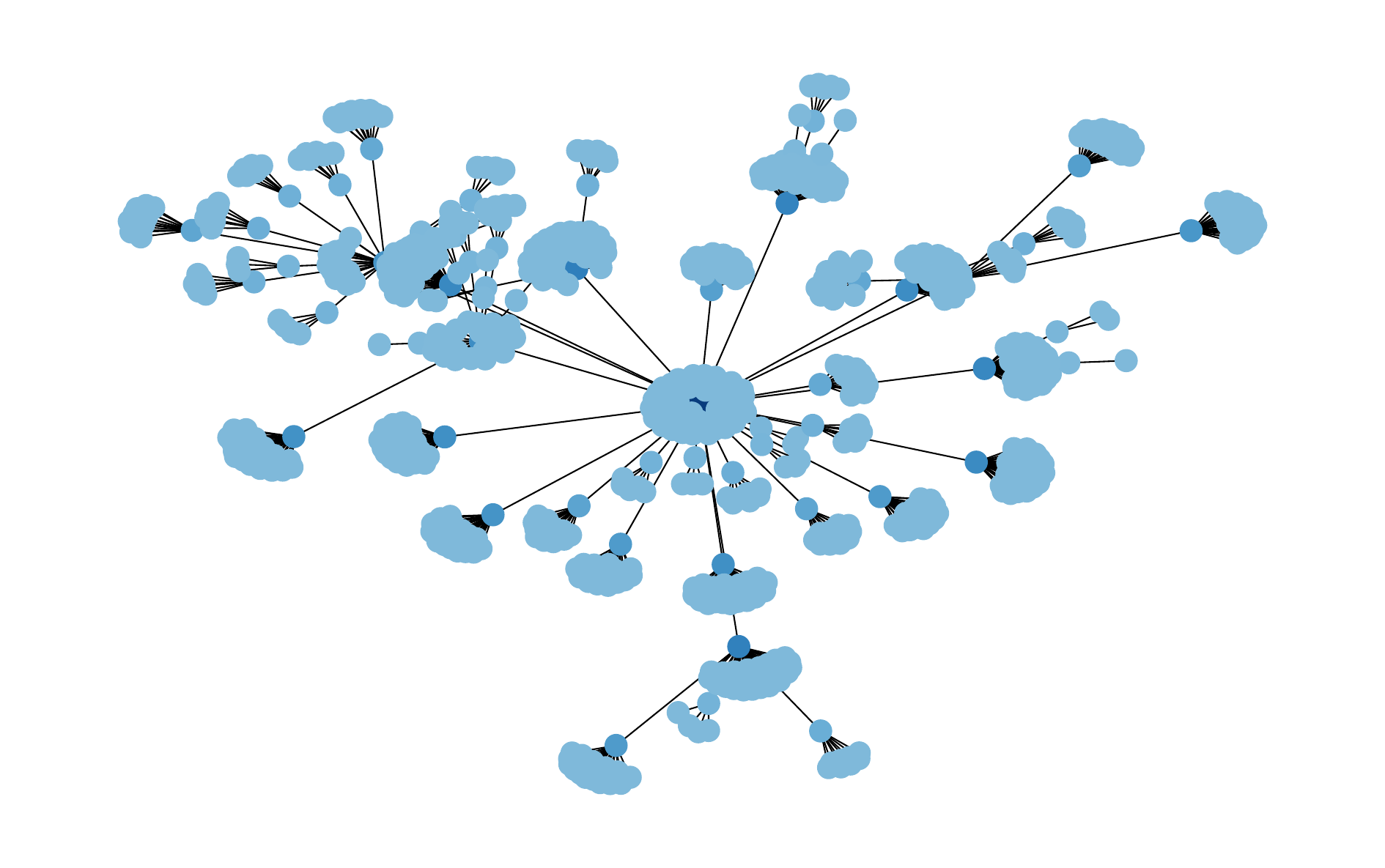}
    \caption{Representation of the Search Space Graph (SSG) for an instance of type $P$}
    \label{fig:solGraph1Instancia}
\end{figure}

\begin{figure}
    \centering
    \includegraphics[width=0.9\linewidth]{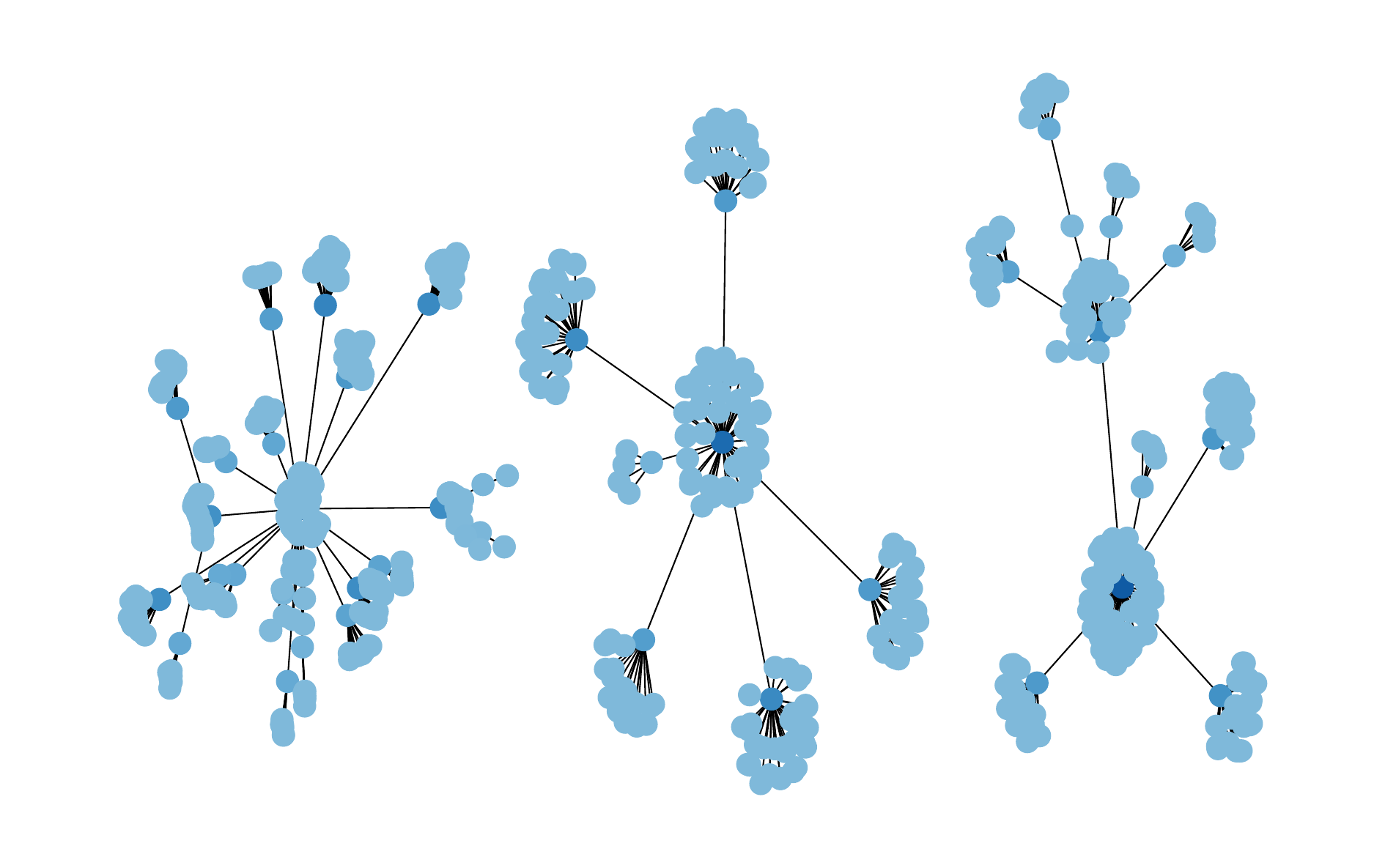}
    \caption{Representation of the Search Space Graph (SSG) for an instance of type $R$}
    \label{fig:solGraph2Instancia}
\end{figure}

This construction brings a method to represent the search space of the instance of the IpMU problem.
Figure \ref{fig:solGraph1Instancia} shows an example of this construction for an instance of $40$ nodes and $100$ edges, where $p=2$ medians must be selected. In this case, the instance belongs to instance type \textit{P}. As can be seen, the graph only represents one connected tree, which means that only one local optimum exists in the search space.
On the other hand, Figure \ref{fig:solGraph2Instancia} represents an instance of type \textit{R} with the very same parameters as the previous one: $n=40$, $m=100$ and $p=2$. As can be seen, in this case, three different connected components exist in the graph, representing the three local optima that can be reached in the search of a solution. It is remarkable that the analysis of larger instances makes the visualization of the different components more difficult.

Although the representation of these graphs with more complex instances could be difficult, we can study the characteristics of SSG without printing the figures. In this case, the study of the disjoint trees gives us an idea of the complexity of the different instances. Contrasting these values with the correlation of cost and time of the instances could give an idea of how this characteristic makes the search space more difficult for the algorithms.

\begin{figure}
    \centering
    \includegraphics[width=1\linewidth]{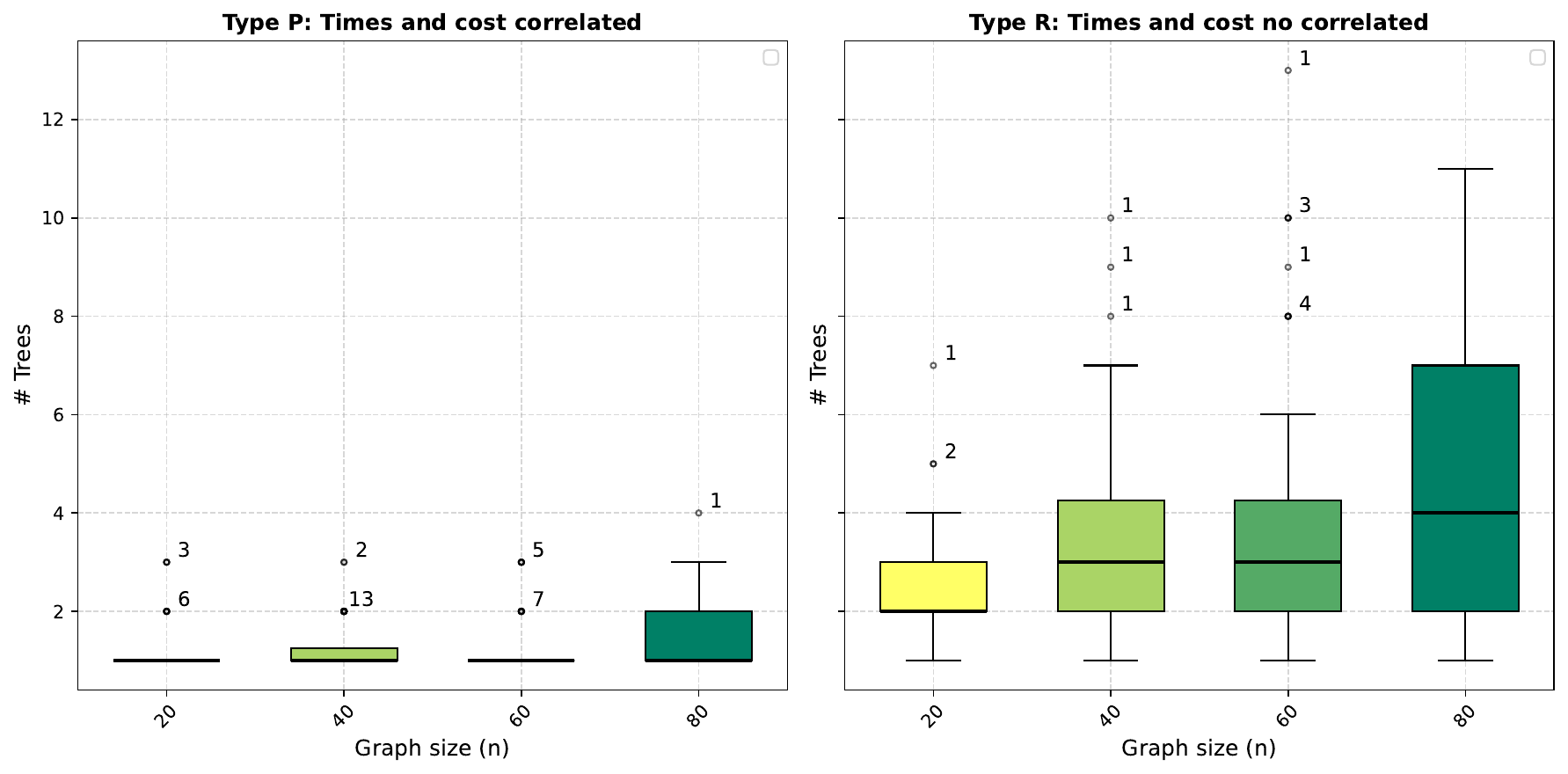}
    \caption{Box plot representing the number of trees in the SSG of each instances divided by size.}
    \label{fig:SSGBoxplot}
\end{figure}

Figure \ref{fig:SSGBoxplot} shows a box plot that shows the evolution of the number of trees (local optimums) in the SSG of each instance. Notice that the set of instances used in this experiment is reduced due to computational limitations. More precisely, the larger instances considered are $n\leq 80$ and $p \leq 3$. The box plot on the left side of the figure shows the results for the set of instances \textit{P}, which present a correlation between costs and time, while the box plot on the right of the figure shows the same experiment for instances \textit{R} with no correlation.

As can be seen in the figure, both types of instances present a growth in the number of trees when the size of the graph increases ($n$). However, it can be seen that the growth is different between both types of instances. The average number of local optimum and the growth is much bigger in the case of non-correlated instances (type \textit{R}), implying a more complex solution space than the other ones. 

As seen in the figure, none of the instances of type \textit{P} presents more than four trees (see size $80$). In contrast, the number of trees is higher for instances of type \textit{R}, reaching up to 13 trees (see size $60$). It should be recalled that this analysis was made in instances with the smallest value of $p$, therefore obtaining a small number of trees.

Making a deeper analysis of the figure, it shows that there is a correlation between the number of local optima existing in the instance problem and the correlation between $C^1$ and $C^2$ values. Both algorithms, metaheuristics and exact approaches, require more computational effort when the number of local optima increases, which explains the difficulty in solving these instances.

A more detailed analysis of the correlated instances shows that at least half of the instances studied in this experiment have just one local (therefore global) optimum. This could explain the execution time spent by the mathematical model to solve these instances since the search space for this instances is convex. On the other hand, it shows that the metaheuristic proposal of this work could reach the optimal value in just one iteration, while in the parameter tuning process performed with \texttt{irace}, a number of 29 iterations were selected for the \textit{GRASP} approach.

Hence, a possible new strategy for the \textit{GRASP} approach could consider the correlation of transportation cost and time (which implies $\mathcal{O}(n)$ complexity) to determine the number of iterations with no improvement. This way, a high correlation will imply just one iteration, saving computation time, and, likely, providing the best solution for the instance.

%This analysis opens a possible new strategy to design the metaheuristic proposal for the IpMU. First, with \texttt{irace} obtain a custom parameter tunning for each set of instances, with and without correlation. Then, to solve a new unknown instance of the IpMU, the algorithm first would start with an analysis of the correlation of both variables, for example evaluating the $R^2$ value, that can be obtained in $\mathcal{O}(n)$ complexity. Then, if both variables are correlated, the probably faster parameter values were established and \textit{GRASP} algorithm is executed; in the other case the same steps could be performed with the other set of parameter values. Although the quality of the results will probably not be affected, the execution time can be considerably reduced.

\section{Conclusions}
\label{sec:conclusions}
The Induced \textit{p}-Median Problem with Upgrades (IpMPU) consists of a combinatorial optimization problem where the best medians have to be selected to provide service to the clients of a network. This problem is highlighted in detaching the concepts of time and cost of transport, using a bi-network graph as input. The problem also introduces some variables related to managing a budget dedicated to upgrading the edges of the graph.

In this work, a metaheuristic algorithm has been proposed to solve the IpMU with performance superior to the state-of-the-art proposals. A two-phase strategy is executed combining a \textit{GRASP} methodology with a greedy algorithm that allows optimal solutions to be obtained in the edge-upgrading subproblem. Our approach is able to obtain the best known result in  {1,226 of the 1,230} studied instances, which corresponds to more than $99.67\%$ of them. We also verified that 1081 of them are optimal solutions of the instances. In contrast, the state-of-the-art model obtains 1108 best results, spending an execution time more than one order of magnitude larger than our \textit{GRASP} proposal in the largest instances.

These results have been also validated with a Bayesian statistical test performed on instances of 100 and 200 nodes. The test shows that our GRASP proposal is able to obtain the best result with a probability of $0.618$, while the state-of-the-art model reaches a probability of $0.380$.

In addition, we have complemented this work by studying the hardness of the problem instances. In particular, we represent the search space of the problem as a directed graph following the trajectory of the local search. The characteristics of the network are contrasted with the characteristics of the instances, showing that the most difficult instances are those with more local optima. This result yields promising results that could speed up the execution of the algorithm and bring about an interesting analysis of the difficulty of each instance of IpMU.

Although the superiority of the GRASP proposal has been demonstrated, there are different ways that could further improve the performance of our algorithmic proposal: the proposed local search method makes an exhaustive search of all the nodes that could be involved in the exchange move, while some heuristic ideas could bring a reduced set of good candidates, further reducing execution time. Following the same idea, enhancing the improving phase of the GRASP proposal with more sophisticated algorithms could improve or even accelerate the proposal.

Finally, we are currently working on using some of the ideas proposed in this paper to other optimization problems of the FLP family. Additionally, the work undertaken to contrast the search space with the characteristics of the instances gives rise to a new analysis of combinatorial optimization problems. Our objective is to establish a more formal methodology that could be applied to different problems.

\section*{Acknowledgements}
This work has been partially supported by the Spanish Ministerio
de Ciencia e Innovaci\'on (MCIN/AEI/10.13039/501100011033) under
grant refs. PID2024-156045NB-I00, TSI-100930-2023-3 and
PID2021-125709OA-C22 and by ERDF A way of making Europe; and Comunidad
Aut\'{o}noma de Madrid with grant refs. PIPF-2024/COM-35036 and
TEC-2024/COM-404.




\bibliographystyle{plainnat}
\bibliography{IpMU}

\end{document}